\documentclass[12pt,a4paper]{article}%
\usepackage[utf8]{inputenc}
\usepackage{hyperref}
\usepackage{amsmath}
\usepackage{amsfonts}
\usepackage{amssymb}
\usepackage{xcolor}
\usepackage{graphicx}%
\providecommand{\U}[1]{\protect\rule{.1in}{.1in}}
\newtheorem{theorem}{Theorem}

\newtheorem{conjecture}[theorem]{Conjecture}
\newtheorem{corollary}[theorem]{Corollary}

\newtheorem{lemma}[theorem]{Lemma}

\newtheorem{problem}[theorem]{Problem}
\newtheorem{proposition}[theorem]{Proposition}

\newenvironment{proof}[1][Proof]{\noindent\textbf{#1.} }{\ \hfill \rule{0.5em}{0.5em}\bigskip}
\graphicspath{{Slike/}}

\begin{document}

\title{On the maximum interval index of a graph}
\author{Martin Knor$^{1}$, Jelena Sedlar$^{2,3}$, Riste \v{S}krekovski$^{3,4}$\\[0.3cm]
{\small $^{1}$ \textit{Slovak University of Technology in Bratislava, Faculty of Civil
Engineering,}}\\[0.05cm]
{\small \textit{Department of Mathematics, Bratislava, Slovakia}}\\[0.1cm]
{\small $^{2}$ \textit{University of Split, Faculty of Civil Engineering, Architecture and
Geodesy, Split, Croatia}}\\[0.1cm]
{\small $^{3}$ \textit{University of Novo mesto, Faculty of Information Studies, Novo mesto, Slovenia}}\\[0.1cm]
{\small $^{4}$ \textit{University of Ljubljana, Faculty of Mathematics and Physics, Ljubljana, Slovenia}}}
\date{}
\maketitle

\begin{abstract}
For a connected graph $G$ on $n$ vertices, the interval index is defined by
$\mathrm{Int}(G)=\sum_{\{u,v\}\subset V(G)}(\left\vert I(u,v)\right\vert -1),$ where $I(u,v)$ denotes the
set of all vertices lying on a shortest $u-v$ path.
It was conjectured that, among all connected graphs on $n$ vertices, the path $P_{n}$ attains the maximum value of the interval index.
In this paper we settle this conjecture, and we prove that the maximum is attained only for $P_{n}$, with the exception of $n=4$, when the cycle $C_{4}$ attains the same value.
We also prove that, for every $n\geq 31$, the maximum interval index among graphs without pendant vertices equals $(n^{3}-n)/6-2(n-4)$ and it is attained by a dumbbell graph, obtained by joining two $4$-cycles by a path.
\end{abstract}



\section{Introduction}

Distance-based graph invariants have been studied for decades, since they capture structural properties of graphs which are relevant, in particular, in mathematical chemistry.
The oldest and most thoroughly investigated such invariant is the Wiener index, introduced by Wiener already in 1947 \cite{Wiener} (see \cite{KST} for a survey of its mathematical aspects).
Recently, Knor, Kovijani\'{c} Vuki\'{c}evi\'{c}, \v{S}krekovski, Do\v{s}ljak, Dimitrov, Popivoda, Vujo\v{s}evi\'{c} and Sedlar \cite{IntervalIndex} introduced the \emph{interval index} of a graph.
This invariant can be viewed as a refinement of the Wiener index, since it takes into account not only the distance between two vertices, but the order of the subgraph formed by the shortest paths joining them.
Several other recent refinements of the Wiener index by (some of) the same authors are in \cite{KSZ,KSY1,KSY2}.

Let $G=(V,E)$ be a connected graph and, for $u,v\in V(G)$, let $d(u,v)$ denote the distance between $u$ and $v$ in $G$.
The \emph{interval} between $u$ and $v$ is
$$
I(u,v)=\{w\in V(G):d(u,w)+d(w,v)=d(u,v)\},
$$
that is, the set of all vertices lying on some shortest $u-v$ path.
This notion goes back to Mulder's monograph \cite{Mulder}.
Nebesk\'{y} \cite{Nebesky} used it to characterize geodetic graphs, that is, graphs in which $\left\vert I(u,v)\right\vert =d(u,v)+1$ for every $u,v$.
The \emph{interval index} of $G$ is
$$
\mathrm{Int}(G)=\sum_{\{u,v\}\subset V(G)}\left( \left\vert I(u,v)\right\vert -1\right) .
$$

If $G$ is a tree, then any two vertices of $G$ are joined by a unique shortest path, so $\mathrm{Int}(G)$ coincides with the Wiener index of $G$.
It is shown in \cite{IntervalIndex} that, among all trees on $n$ vertices, the path $P_{n}$ maximizes and the star $S_{n}$ minimizes the interval index.
More generally, a graph with a cut vertex can be decomposed according to a formula which allows the interval index to be computed from its blocks.
Further, the extremal values of the interval index were determined for unicyclic graphs of given girth,
and $\mathrm{Int}(G)$ was computed for hypercubes, for complete bipartite graphs and for ladder graphs.
Since the complete graph minimizes the interval index among all connected graphs on $n$ vertices, while the maximum turned out to be far less accessible, the authors of \cite{IntervalIndex} posed the following conjecture.

\begin{conjecture}
\label{Con_path}
In the class of connected graphs on $n$ vertices the path $P_{n}$ attains the maximum value
of interval index.
\end{conjecture}

The analogous question for graphs without pendant vertices is more delicate.
Extremal problems for the Wiener index under a minimum-degree constraint were studied by Alochukwu and Dankelmann \cite{DA}, and for $2$-connected graphs by Bessy, Dross, Knor and \v{S}krekovski \cite{BDKS}.
The problem which graphs of a given minimum degree attain the maximum Wiener index seems to be open in general, see Problem~4.1 in \cite{KST}.
Among $2$-connected graphs on $n$ vertices the maximum Wiener index is attained by the cycle $C_{n}$, see \cite{BDKS} and the references therein, but this is no longer true for graphs with minimum degree at least $2$.
In \cite{BDHKS} it is determined which graph maximizes the Wiener index in the class of connected graphs with $p$ blocks.
From these results it follows that in the class of graphs without pendant vertices the maximum is attained by two triangles joined by a path.
To see this, it is sufficient to realize that connected graphs which have at least $n-2$ blocks have a pendant edge (they consist of a path with one triangle if $p=n-2$, and the path is the unique graph with $p=n-1$), and to compare the Wiener index of a cycle with that of two triangles joined by a path.
However, the graph consisting of two triangles joined by a path is not extremal for the interval index.
Namely, let $x$ and $y$ be the two vertices of one of its triangles which do not lie on the path.
Every other vertex $z$ is at the same distance from $x$ and $y$, and $d(x,y)=1$, so none of the vertices $x,y,z$ lies on a shortest path between the other two.
The same holds for the other triangle, and no triple of vertices contains both such pairs.
Hence, by Lemma~\ref{Lemma_identity} below, the interval index of this graph is at most $(n^{3}-n)/6-2(n-2)$.
For $n\geq 8$ this is less than the value $(n^{3}-n)/6-2(n-4)$ attained by the dumbbell graph $D_{n}$ from Proposition~\ref{Prop_dumbbell}.
So the authors of \cite{IntervalIndex} posed the following problem.

\begin{problem}
\label{Prob_pendant}
Characterize connected graphs on $n$ vertices without pendant vertices which maximize the interval index.
\end{problem}

In this paper we settle Conjecture \ref{Con_path}.
We show that the path $P_n$ is the unique graph attaining the maximum, with the exception of $n=4$, when it is attained also by $C_4$.
Moreover, for $n\geq 31$ we determine the maximum interval index in the class of connected graphs without pendant vertices: it equals $(n^{3}-n)/6-2(n-4)$.
It is attained by the dumbbell consisting of two $4$-cycles joined by a path.

\section{General graphs}

Let $G$ be a connected graph on $n$ vertices with vertex set $V(G)$ and distance function $d$.
Recall that $I(u,v)=\{w\in V(G):d(u,w)+d(w,v)=d(u,v)\}$.
For three pairwise distinct vertices $x,y,z$ of $G$, we say that $y$ \emph{lies between} $x$ and $z$ if $d(x,y)+d(y,z)=d(x,z)$, that
is, if $y\in I(x,z)\setminus \{x,z\}$.
We call an unordered triple $\{x,y,z\}$ of distinct vertices \emph{collinear} if at least one of its three vertices lies between the other two, and we denote by $\mathrm{col}(G)$ the number of collinear triples in $G$.

Notice that a triple may a priori be collinear in three different ways, since each of its three
vertices could be the one lying between the other two.
We show that this is not the case and a triple can be collinear in a unique way.

\begin{lemma}
\label{Lemma_betweenness}
Among three pairwise distinct vertices of a connected graph, at most
one lies between the other two.
\end{lemma}

\begin{proof}
Assume to the contrary that two of the three vertices lie between the other two. Without loss of generality, we may assume that $y$ lies between $x$ and $z$, and $x$ lies between $y$ and $z$.
This means that
$$
d(x,y)+d(y,z)=d(x,z)\quad \text{and}\quad d(y,x)+d(x,z)=d(y,z).
$$
Substituting the first equality into the second one gives $2d(x,y)=0$.
So $x=y$, which is a contradiction with $x,y,z$ being pairwise distinct.
\end{proof}

The uniqueness of the middle vertex is what makes the interval index accessible to counting.

\begin{lemma}
\label{Lemma_identity}
For every connected graph $G$ on $n$ vertices it holds that
$$
\mathrm{Int}(G)=\binom{n}{2}+\mathrm{col}(G).
$$
\end{lemma}

\begin{proof}
By Lemma~{\ref{Lemma_betweenness}}, in a collinear triple the middle vertex (and hence also the two non-middle vertices) is determined uniquely.
So split the collinear triples according to the non-middle vertices.
For non-middle vertices $u$ and $v$ the triple $\{u,v,x\}$ is collinear if and only if $x\in I(u,v)\setminus\{u,v\}$.
So there are $\left\vert I(u,v)\right\vert -2$ collinear triples having $u$ and $v$ as non-middle vertices.
Consequently
$$
\mathrm{Int}(G)=\sum_{\{u,v\}\subset V(G)}(\left\vert I(u,v)\right\vert -1)=\binom{n}{2}+\mathrm{col}(G)
$$
as required.
\end{proof}

We are now in a position to prove Conjecture \ref{Con_path}.

\begin{theorem}
\label{Tm_path}For every connected graph $G$ on $n\geq 1$ vertices it holds that
$$
\mathrm{Int}(G)\leq \binom{n}{2}+\binom{n}{3}=\frac{n^{3}-n}{6}=\mathrm{Int}(P_{n}),
$$
with equality if and only if $G\cong P_{n},$ or $n=4$ and $G\cong C_{4}.$
\end{theorem}

\begin{proof}
The number $\mathrm{col}(G)$ counts a subfamily of the $\binom{n}{3}$ triples of $G$.
So Lemma \ref{Lemma_identity} gives $\mathrm{Int}(G)=\binom{n}{2}+\mathrm{col}(G)\leq \binom{n}{2}+\binom{n}{3}$.
Notice that in the path $P_{n}$ every triple is collinear, since for any three vertices of a path the middle one lies between the other two.
Hence $\mathrm{col}(P_{n})=\binom{n}{3}$ and
$\mathrm{Int}(P_{n})=\binom{n}{2}+\binom{n}{3}$.
Consequently, $P_n$ maximizes the interval index.

Now we prove the uniqueness.
If a graph $G$ on $n$ vertices attains $\mathrm{Int}(P_n)$, then every triple
of $G$ must be collinear.
For $n\leq 2$ there is nothing to prove, so assume that $n\geq 3$.
Let $u,v$ be a diametral pair of $G$, that is, a pair of vertices at distance $\mathrm{diam}(G)$, and denote $D=d(u,v)=\mathrm{diam}(G)$. For
$0\leq i\leq D$ let $L_{i}=\{x\in V(G):d(u,x)=i\}$.
We establish the theorem through several claims.

\medskip

\noindent\textbf{Claim A.}
\emph{If $w\in L_{i}$ then $d(w,v)=D-i$.}

\smallskip

\noindent Let $w\not \in \{u,v\}$. The triple $\{u,v,w\}$ is collinear, so one of its vertices
lies between the other two. Assume that $u$ lies between $w$ and $v$.
Then $d(w,v)=d(w,u)+d(u,v)=d(w,u)+D>D$, which is a contradiction with $D=\mathrm{diam}(G)$.
The case when $v$ lies between $w$ and $u$ is symmetric.
Hence $w$ lies between $u$ and $v$, that is, $d(u,w)+d(w,v)=D$.

\medskip

\noindent\textbf{Claim B.}
\emph{If $1\leq i\leq D-1$ and $w\not =w^{\prime}$ both belong to $L_{i}$, then $i=D/2$ and $d(w,w^{\prime})=2i$.}

\smallskip

\noindent Consider the triple $\{u,w,w^{\prime}\}$.
Here $d(u,w)=d(u,w^{\prime})=i$.
So neither $w$ nor $w^{\prime}$ can lie between $u$ and the other vertex, as this would force $d(w,w^{\prime})=0$.
Hence $u$ lies between $w$ and $w^{\prime}$, which gives $d(w,w^{\prime})=2i$.
By Claim A we have $d(w,v)=d(w^{\prime},v)=D-i$. So the same argument applied to the triple $\{v,w,w^{\prime}\}$ gives $d(w,w^{\prime})=2(D-i)$.
Therefore $2i=2(D-i)$, and hence $i=D/2$.
\medskip

By Claim B, every level $L_{i}$ with $1\leq i\leq D-1$ and $i\not =D/2$ is a singleton.
By Claim A, $L_{D}=\{v\}$, and of course $L_{0}=\{u\}$.
Consequently, if $D$ is odd then all levels are singletons, i.e., $G$ is the path.

\medskip

\noindent\textbf{Claim C.}
\emph{If some level $L_{i}$ contains at least two vertices, then $D=2$, $n=4$ and $G\cong C_{4}$.}

\smallskip

\noindent By Claim B this can happen only for $i=D/2\geq 1$.
Denote $m=\left\vert L_{i}\right\vert\geq 2$, and assume first that $m\geq 3$.
Let $w_{1},w_{2},w_{3}\in L_{i}$ be pairwise distinct.
By Claim B they are pairwise at distance $D$.
So none of them can lie between the other two, as $D+D=D$ would force $D=0$.
Hence the triple $\{w_{1},w_{2},w_{3}\}$ is not collinear,
which is a contradiction with our assumption.
So $m=2$, and we denote $L_{i}=\{w,w^{\prime}\}$.
Since $i-1\not =D/2$, Claim B makes $L_{i-1}=\{a\}$ a singleton.
Both $w$ and $w^{\prime}$ have a neighbour in $L_{i-1}$, necessarily $a$.
Hence $d(w,w^{\prime})\leq 2$. Claim B gives $d(w,w^{\prime})=2i$, so $i=1$. 
Therefore $D=2$ and $n=4$.
Since $d(u,v)=2$ and $d(w,w^{\prime})=2$, $G\cong C_{4}$.

To summarize, either all levels are singletons and $G\cong P_{n}$, or Claim C applies and $G\cong C_{4}$.
Conversely, $\mathrm{Int}(C_{4})=4\cdot 1+2\cdot 3=10=\mathrm{Int}(P_{4})$.
\end{proof}

\section{Graphs without pendant vertices}

Recall that Problem \ref{Prob_pendant} asks for the graphs on $n$ vertices without
pendant vertices, that is, with $\delta(G)\geq 2,$ which maximize the interval index.

In this section we determine the maximum interval index in Problem~{\ref{Prob_pendant}} for $n\geq 31$.
A triple which is not collinear will be called \emph{non-collinear}.
Non-collinear triples satisfy the strict triangle inequality on all three of their pairwise distances.
For a connected graph $G$ on $n$ vertices we denote the number of its non-collinear triples by $\mathrm{noncol}(G)$.
Then by Lemma \ref{Lemma_identity} we have
$$
\mathrm{Int}(G)=\binom{n}{2}+\binom{n}{3}-\mathrm{noncol}(G).
$$
Hence, maximizing the interval index over graphs without pendant vertices is the same as minimizing the number of non-collinear triples over the same class.
Notice that the number of non-collinear triples is what a graph pays for not being a path.

For $n\geq 8$, the \emph{dumbbell} $D_{n}$ is the graph obtained by joining two disjoint copies of $C_{4}$ by a path of length $n-7$, that is with $n-8$ internal vertices.
We show that $D_n$ is an extremal graph in the class of graphs without pendant vertices.

\begin{proposition}
\label{Prop_dumbbell}
For every $n\geq 8$ it holds that $\mathrm{noncol}(D_{n})=2(n-4)$.
\end{proposition}

\begin{proof}
Denote the two copies of $C_{4}$ by $a^{\prime},w_{1},a,w_{1}^{\prime}$ and $b^{\prime},w_{2},b,w_{2}^{\prime}$.
Here $a$ and $b$ are the vertices at which the connecting path is attached, while $a^{\prime}$ and $b^{\prime}$ are the vertices antipodal to them.

Take any $z\not \in \{w_{1},w_{1}^{\prime}\}$.
Then $d(z,w_{1})=d(z,w_{1}^{\prime})$, and $d(w_{1},w_{1}^{\prime})=2$.
So the triple $\{w_{1},w_{1}^{\prime},z\}$ is collinear if and only if $z$ lies between $w_{1}$ and $w_{1}^{\prime}$.
Hence $z\in \{a,a^{\prime}\}$.
Consequently, exactly $n-4$ triples of this form are non-collinear.
The pair $\{w_{2},w_{2}^{\prime}\}$ gives another $n-4$ non-collinear triples $\{w_2,w_2^{\prime},z\}$ disjoint from the previous ones.

Every remaining triple avoids at least one of $w_{1},w_{1}^{\prime}$ and at least one of $w_{2},w_{2}^{\prime}$.
Delete one vertex from each pair.
What is left is an isometric path in $D_{n}$ containing the triple.
All triples of such a path are collinear, so $\mathrm{noncol}(D_{n})=2(n-4)$.
\end{proof}

In order to show that no graph without pendant vertices beats $D_{n}$, we distinguish whether $G$ has a cut vertex or not.
Let us first consider the case in which $G$ has a cut vertex.
Since $G$ has no pendant vertices, no pendant block of $G$ is a bridge.
So every leaf of the block tree of $G$ is a $2$-connected block $B$ on $k\geq 3$ vertices, and we denote by $a$ the unique cut vertex of $G$ lying in $B$.
Notice that distances inside a block agree with distances in $G$, so we may measure them in $B$ alone.

A pair of vertices of $B$ other than $a$ is \emph{bad} at $a$ if neither of the two vertices lies on a shortest path from $a$ to the other one.
That is, a pair $\{x,y\}\subseteq V(B)\setminus\{a\}$ is bad if $d(x,y)\not =\left\vert d(a,x)-d(a,y)\right\vert$.
We denote by $\beta(B,a)$ the number of bad pairs of $B$ at $a$.

\begin{lemma}
\label{Lemma_leafblock}
Let $B$ be a leaf block of $G$ on $k$ vertices with cut vertex $a$.
Then every bad pair of $B$ at $a$ forms a non-collinear triple with each of the $n-k$ vertices of $G$ outside $B$.
So there are at least
$$
\beta(B,a)(n-k)+\mathrm{noncol}(B)
$$
non-collinear triples having at least two vertices in $V(B)\setminus\{a\}$.
\end{lemma}

\begin{proof}
Let $\{x,y,z\}$ be a triple such that $x,y\in V(B)\setminus\{a\}$ and $z\in V(G)\setminus V(B)$.
Every path from $x$ or $y$ to $z$ passes through $a$, so $d(x,z)=d(x,a)+d(a,z)$ and $d(y,z)=d(y,a)+d(a,z)$.
Assume that $\{x,y,z\}$ is collinear.
Then $z$ cannot be the middle vertex since $d(x,z)+d(y,z)>d(x,a)+d(y,a)\geq d(x,y)$.
But if $x$ is the middle vertex then
$$
d(z,a)+d(a,y)=d(z,y)=d(z,x)+d(x,y)=d(z,a)+d(a,x)+d(x,y)
$$
so $d(a,y)=d(a,x)+d(x,y)$ and hence $\{x,y\}$ is not a bad pair at $a$.
The case when $y$ is the middle vertex is symmetric.
Thus every bad pair forms a non-collinear triple with each vertex outside $B$.

The triples inside $B$ are counted by $\mathrm{noncol}(B)$.
They remain non-collinear in $G$, since distances inside $B$ agree with those in $G$.
They contain at most one vertex $a$, and they are different from the triples above, which contain a vertex outside $B$.
\end{proof}

\begin{lemma}
\label{Lemma_bounds}
Let $B$ be a $2$-connected graph on $k\geq 3$ vertices and $a$ a vertex of $B$.
Then
$$
\beta(B,a)+\mathrm{noncol}(B)\geq k-3.
$$
\end{lemma}

\begin{proof}
Denote by $e$ the eccentricity of $a$ in $B$ and by $D$ the diameter of $B$.
We establish the lemma through two claims.
The first one bounds the number of bad pairs at $a$, and the second one the number of non-collinear triples of $B$.

\medskip

\noindent\textbf{Claim A.}
\emph{We have $\beta(B,a)\geq k-1-e$.}

\smallskip

\noindent Let $L_{i}=\{x\in V(B):d(a,x)=i\}$ for $0\leq i\leq e$, and put $\ell_{i}=\left\vert L_{i}\right\vert$.
Take any two vertices $x,y$ of the same level $L_{i}$.
Then $d(a,x)=d(a,y)$ and $d(x,y)\geq 1$, so the pair $\{x,y\}$ is bad.
Hence
$$
\beta(B,a)\geq \sum_{i=1}^{e}\binom{\ell_{i}}{2}\geq \sum_{i=1}^{e}(\ell_{i}-1)=k-1-e.
$$

\medskip

\noindent\textbf{Claim B.} \emph{We have $\mathrm{noncol}(B)\geq D-2$.}

\smallskip

\noindent Let $u,v$ be a diametral pair in $B$, and let $M_{i}=\{x\in V(B):d(u,x)=i\}$ for $0\leq i\leq D$.
Then $\left\vert M_{i}\right\vert \geq 2$ for $1\leq i\leq D-1$, since otherwise the unique vertex of $M_{i}$ would be a cut-vertex separating $u$ from $v$.
Fix an $i$ such that $i\not =D/2$, and pick two distinct vertices $w,w^{\prime}\in M_{i}$.
We show that one of the triples $\{w,w^{\prime},u\}$ and $\{w,w^{\prime},v\}$ is not collinear.

We have $d(u,w)=d(u,w^{\prime})=i$, so neither $w$ nor $w^{\prime}$ can lie between the other two vertices of $\{u,w,w^{\prime}\}$.
Hence this triple is collinear only if $u$ lies between $w$ and $w^{\prime}$, that is, only if $d(w,w^{\prime})=2i$.
For $2i>D$ we have $d(w,w^{\prime})\leq D<2i$, so the triple $\{u,w,w^{\prime}\}$ is non-collinear.

Assume therefore that $2i<D$, and assume that both $\{u,w,w^{\prime}\}$ and $\{v,w,w^{\prime}\}$ are collinear.
The first one gives $d(w,w^{\prime})=2i$.
From $D=d(u,v)\leq d(u,w)+d(w,v)=i+d(v,w)$ we get $d(v,w)\geq D-i$, and likewise $d(v,w^{\prime})\geq D-i$.
Assume that $w$ is between $v$ and $w^{\prime}$.
Then $d(v,w^{\prime})=d(v,w)+d(w,w^{\prime})\geq (D-i)+2i>D$, a contradiction.
The case in which $w^{\prime}$ is between $v$ and $w$ is symmetric.
Hence $v$ is between $w$ and $w^{\prime}$, which gives $2(D-i)\leq d(v,w)+d(v,w^{\prime})=d(w,w^{\prime})=2i$.
So $2i\geq D$, a contradiction.
Therefore at least one of the two triples is non-collinear.

In both cases we obtain a non-collinear triple.
It consists of the pair $\{w,w^{\prime}\}\subseteq M_{i}$ together with one of $u,v$.
Since distinct values of $i$ give distinct pairs and distinct triples, we get
$$
\mathrm{noncol}(B)\geq \left\vert \{1\leq i\leq D-1:i\not =D/2\}\right\vert \geq (D-1)-1=D-2.
$$

\medskip

\noindent The eccentricity of a vertex never exceeds the diameter, so $D\geq e$. Claims A and B now give $\beta(B,a)+\mathrm{noncol}(B)\geq (k-1-e)+(D-2)\geq k-3$.
\end{proof}

\begin{theorem}
\label{Tm_cutvertex}
Every connected graph on $n$ vertices without pendant vertices which has a cut vertex has at least $2(n-4)$ non-collinear triples.
\end{theorem}

\begin{proof}
The block tree of $G$ has at least two leaves. Let $B_{1},B_{2}$ be two distinct leaf blocks, with cut vertices $a_{1},a_{2}$ and orders $k_{1},k_{2}$.
Two distinct blocks of $G$ meet in at most one vertex, and such a common vertex would be a cut vertex lying in both leaf blocks, so it would be equal to $a_{1}$ and $a_{2}$.
Hence the sets $V(B_{1})\setminus \{a_{1}\}$ and $V(B_{2})\setminus \{a_{2}\}$ are disjoint.
A triple counted by Lemma \ref{Lemma_leafblock} for both blocks would need at least two vertices in each of these two sets.
Hence the two collections of triples are disjoint, which gives
$$
\mathrm{noncol}(G)\geq \sum_{i=1,2}\left[ \beta(B_{i},a_{i})(n-k_{i})+\mathrm{noncol}(B_{i})\right].
$$

It suffices to show that $f:=\beta(n-k)+\mathrm{noncol}(B)\geq n-4$ for a single leaf block $B$ of order $k$ with cut vertex $a$, where we abbreviate $\beta=\beta(B,a)$.
There is at least one vertex of $G$ outside $B$, and so $n-k\geq 1$.
Notice that $\beta\geq 1$.
Indeed, $B$ is $2$-connected on $k\geq 3$ vertices, so $a$ has two neighbours $x,y$ in $B$.
Then $d(a,x)=d(a,y)=1$ and $d(x,y)\geq 1$, so $\{x,y\}$ is a bad pair.

Therefore
\begin{align*}
f-(n-4)  &  =\beta(n-k)+\mathrm{noncol}(B)-(n-4) = (\beta-1)(n-k)+\mathrm{noncol}(B)-(k-4) \\
&  \geq(\beta-1)+\mathrm{noncol}(B)-(k-4) = \left[ \beta+\mathrm{noncol}(B)\right] -(k-3)\geq 0,
\end{align*}
where in the first inequality we used $n-k\geq 1$ and $\beta-1\geq 0$, and in the last one Lemma \ref{Lemma_bounds}.
Hence $f\geq n-4$, which finishes the proof.
\end{proof}

By Proposition \ref{Prop_dumbbell} the bound of the previous theorem is attained, so among the graphs without pendant vertices which have a cut vertex the dumbbell is extremal for every $n\geq 8$.
It remains to treat the $2$-connected case.
First we solve the case when the graph has diameter at most $3$.

\begin{theorem}
\label{Tm_diam3}
Every connected graph $G$ on $n\geq 8$ vertices of diameter at most $3$ has at least $2(n-4)$ non-collinear triples.
\end{theorem}

\begin{proof}
Denote by $d_{x}$ the degree of a vertex $x$ of $G$.
We say that a set $S\subseteq V(G)$ \emph{spans} the edges of $G$ with both end-vertices in $S$.
We first establish two claims. The first one tells which triples can be collinear at all, and the second one counts them.

\medskip

\noindent\textbf{Claim A.} \emph{A triple of $G$ which spans no edge, or which spans all three of them, is not collinear.}

\smallskip

\noindent A triple spanning three edges is a triangle.
All its three distances equal $1$, and $1+1>1$, so it is not collinear.
Take now a triple spanning no edge.
Its three distances are at least $2$, and at most $3$ because $\mathrm{diam}(G)\leq 3$.
Here $3<2+2$, so this triple is not collinear either.

\medskip

\noindent\textbf{Claim B.} \emph{The number of triples of $G$ spanning exactly one or exactly two edges equals one half of the sum $\sum_{x\in V(G)}d_{x}(n-1-d_{x}).$}\smallskip

\noindent Call a vertex of a triple \emph{mixed} if inside that triple it is joined to exactly one of the other two vertices.
A triple spanning exactly one edge has exactly two mixed vertices, namely the ends of that edge.
A triple spanning exactly two edges also has exactly two mixed vertices, namely the ends of the missing pair.
Triangles and edgeless triples have none. So a triple is counted by the claim if and only if it carries mixed vertices, and then it carries precisely two of them.
It remains to count the pairs $(x,\{y,z\})$ in which $x$ is mixed, that is, in which $x$ is adjacent to $y$ and not to $z$.
For a fixed $x$ there are $d_{x}(n-1-d_{x})$ such pairs.
We sum over $x$ and divide by two, since each considered triple contains exactly two mixed vertices.

\medskip

\noindent Claims A and B together give
$$
\mathrm{noncol}(G)\geq \binom{n}{3}-\frac{1}{2}\sum_{x\in V(G)}d_{x}(n-1-d_{x}),
$$
so it remains to bound the sum from above.

Consider the quadratic function $q(x)=x(n-1-x)$.
We have $q(x)\leq (n-1)^{2}/4$, and for even $n$ and integer $x$ we get $q(x)\leq n(n-2)/4$.
This yields $\mathrm{noncol}(G)\geq n(n-1)(n-5)/24$ for odd $n$, and $\mathrm{noncol}(G)\geq n(n-2)(n-4)/24$ for even $n$.
Both bounds are at least $2(n-4)$ for $n\geq 8$.
Namely, for even $n$ this is equivalent to $n(n-2)\geq 48$, which holds for $n\geq 8$.
For odd $n\geq 9$ we have $n(n-1)(n-5)\geq 72(n-5)\geq 48(n-4)$.
\end{proof}

Notice that the proof gives for graphs of diameter at most $3$ a bound which is cubic in $n$, hence much stronger than the linear one we are after.
Thus it is only the graphs of large diameter which may come close to the dumbbell.

Throughout the rest of this section $G$ is a $2$-connected graph on $n$ vertices and $D=\mathrm{diam}(G)\geq 4$.
The vertices at distance $i$ from a vertex $u$ form the \emph{$i$-th level} of $u$,
\[
M_{i}(u)=\{x\in V(G):d(u,x)=i\},\qquad m_{i}(u)=\left\vert M_{i}(u)\right\vert .
\]
A pair $(u,v)$ with $d(u,v)=D$ is \emph{diametral}.
We fix one such pair and abbreviate $m_{i}=m_{i}(u)$.
We write $m^{\ast}=m_{D/2}(u)$ if $D$ is even, $m^{\ast}=0$ otherwise.
A pair of distinct vertices $w,w^{\prime}$ of a level $M_{i}(u)$ is \emph{$u$-good} if $u$ does not lie between them, that is, if $d(w,w^{\prime})\not =2i$.

In the following statement we reformulate the notion of collinearity, in which one vertex of the triple is singled out as an \emph{apex}.
It says that a triple is collinear precisely when the distances from the apex to the other two vertices either add up to, or differ by, the distance
between these other two vertices.
Part $(i)$ of the next statement is referred to as the \emph{apex criterion}.

\begin{lemma}
\label{Lemma_basic}
Let $G$ be a $2$-connected graph of diameter $D$, let $(u,v)$ be a diametral pair and $m_i=\left\vert M_i(u)\right\vert$, and let $z,w,w^{\prime}$ be three pairwise distinct vertices of $G.$ The following statements hold:

\begin{enumerate}
\item[\textrm{(i)}] denoting $a=d(z,w),$ $b=d(z,w^{\prime})$ and $r=d(w,w^{\prime})$, the triple $\{z,w,w^{\prime}\}$ is collinear if and only if $a+b=r$ or $\left\vert a-b\right\vert =r$;

\item[\textrm{(ii)}] if $w$ and $w^{\prime}$ belong to a level $M_{i}(u)$ with $2i>D,$ then the pair $\{w,w^{\prime}\}$ is $u$-good and the triple $\{u,w,w^{\prime}\}$ is not collinear;

\item[\textrm{(iii)}] we have $m_{i}\geq 2$ for $1\leq i\leq D-1$.
\end{enumerate}
\end{lemma}

\begin{proof}
$(i)$ The three betweenness relations read: $z$ lies between $w$ and $w^{\prime}$ if and only if $a+b=r$; $w$ lies between $z$ and $w^{\prime}$ if and only if $a+r=b$; and $w^{\prime}$ lies between $z$ and $w$ if and only if $b+r=a$.
The last two together say $\left\vert a-b\right\vert=r$.

$(ii)$ Here $a=b=i$.
So by $(i)$ the triple is collinear if and only if $d(w,w^{\prime})=2i$ or $d(w,w^{\prime})=0$.
The first fails because $d(w,w^{\prime})\leq D<2i$, the second because $w\not =w^{\prime}$.

$(iii)$ If $m_i=1$, where $1\leq i\leq D-1$, then the single vertex in $M_i(u)$ separates $u$ from $v$, since an edge changes the distance from $u$ by at most one, so every $u$--$v$ path meets every level $M_i(u)$, $1\leq i\leq D-1$.
This contradicts the $2$-connectivity of $G$.
\end{proof}

As in the case of graphs with a cut vertex, we bound the number of non-collinear triples.
Let $M_0(u)=\{u\}$.
We consider the following four types of triples:

\begin{enumerate}
\item[$\tau_1$:]
Triples where all three vertices are from one level, namely the middle one $M_{D/2}(u)$, where $D$ is even.

\item[$\tau_2$:]
Triples $\{x,w,w^{\prime}\}$, where $w,w^{\prime}\in M_{i}(u)$ for some $1\leq i\leq D$ and $x\in \{u,v\}\setminus M_{i}(u)$.
In particular, if $i=D$ then $x=u$.

\item[$\tau_3$:]
Triples containing two vertices $w,w^{\prime}$ of a common level and a vertex $z$ from a different level, $z\not \in \{u,v\}$.

\item[$\tau_4$:]
Triples of vertices from three distinct levels $M_i(u)$, $0\leq i\leq D$.
\end{enumerate}

Notice that a triple of type $\tau_4$ has no two vertices in a common level, while a triple of type $\tau_1$ has all three.
The triples of types $\tau_2$ and $\tau_3$ differ in whether the third vertex belongs to $\{u,v\}$.
Hence the types are pairwise disjoint, and we count only these triples; for instance, three vertices of a level other than the middle one have no type.
First we bound the number of non-collinear triples of type $\tau_1$ and $\tau_2$, then the number of those of type $\tau_3$ and $\tau_4$.

For estimating the number of triples in the middle level we need the following statement.
It is formulated for $\mu\geq 0$, since the corresponding graph can be small, even empty.
For a graph $J$, by $t(J)$ we denote the number of triangles in $J$.

\begin{lemma}
\label{Lemma_triangles}
Every graph $J$ on $\mu\geq 0$ vertices satisfies
$$
\binom{\mu}{2}-2\left\vert E(J)\right\vert +t(J)\geq -\frac{\mu+1}{2}.
$$

\end{lemma}

\begin{proof}
Write $e=\left\vert E(J)\right\vert $ and assume first that $e\leq \mu^{2}/4$. Then the first two terms alone give $\binom{\mu}{2}-2e\geq \mu(\mu-1)/2-\mu^{2}/2=-\mu/2$.
Observe that this covers the cases $\mu\leq 2$.

Assume next that $e>\mu^{2}/4$.
The supersaturation bound of Moon and Moser \cite{MoonMoser} gives $t(J)\geq e(4e-\mu^{2})/(3\mu)$.
So the left-hand side is at least $g(e)=\binom{\mu}{2}-2e+e(4e-\mu^{2})/(3\mu)$, and $g$ is convex with $g^{\prime}(e)=-2+(8e-\mu^{2})/(3\mu)$.
For $\mu\geq 6$ we get $g^{\prime}(\mu^{2}/4)=\mu/3-2\geq 0$.
Hence $g$ is non-decreasing on $[\mu^{2}/4,\infty)$, which gives $g(e)\geq g(\mu^{2}/4)=-\mu/2$.
For $\mu\leq 5$ the constraint $\mu^{2}/4<e\leq \binom{\mu}{2}$ leaves only finitely many pairs $(\mu,e)$, namely $(3,3)$, $(4,5)$, $(4,6)$, $(5,7)$, $(5,8)$, $(5,9)$ and $(5,10)$.
For them the values of $g(e)$ are $-2$ for $\mu=3$, then $-7/3$ and $-2$ for $\mu=4$, and $-13/5$, $-34/15$, $-7/5$ and $0$ for $\mu=5$.
In each case this is at least $-(\mu+1)/2$, which equals $-2$, $-5/2$ and $-3$ for $\mu=3,4,5$, respectively.
\end{proof}

The middle level requires some more notation.
If $D$ is even, let $H$ be a graph on the vertex set $M_{D/2}(u)$ in which two vertices are joined whenever they are at distance exactly $D$ in $G$.
Let $M=\{w\in V(H):d(v,w)=D/2\}$ be the set of those vertices of $M_{D/2}(u)$ which are as far from $v$ as they are from $u,$ and let $J=H[M]$ be the subgraph of $H$ induced by $M$.

\begin{lemma}
\label{Lemma_levels}Let $G$ be a $2$-connected graph of diameter $D\geq 4$, and let $(u,v)$ be a diametral pair of $G$.
Then the number of non-collinear triples of type $\tau_1$ and $\tau_2$ is at least
$$
\sum_{i=1}^D\binom{m_i}{2}-\frac{m^{\ast}+1}2.
$$
\end{lemma}

\begin{proof}
We consider a more subtle partition of non-collinear triples of types $\tau_1$ and $\tau_2$.
In each case we bound their number and at the end we sum these bounds.

\medskip

\noindent {\bf Case 1}:
{\em Triples containing $w,w^{\prime}\in M_i(u)$, where $1\leq i\leq D-1$ and $i\not =D/2$.}

\smallskip

Suppose that both triples $\{u,w,w^{\prime}\}$ and $\{v,w,w^{\prime}\}$ are collinear.
Since $d(u,w)=d(u,w^{\prime})=i$, the apex criterion at $u$ forces $d(w,w^{\prime})=2i$.
Denote $j=d(v,w)$ and $j^{\prime}=d(v,w^{\prime})$.
The apex $v$ gives $d(w,w^{\prime})=2i=j+j^{\prime}$ or $2i=\left\vert j-j^{\prime}\right\vert $.
From $D=d(u,v)\leq d(u,w)+d(w,v)$ we get $j\geq D-i$, and likewise $j^{\prime}\geq D-i$. Of
course, $j,j^{\prime}\leq D$.
Assume that $2i=j+j^{\prime}$.
Then $2i=j+j^{\prime}\geq 2(D-i)$ gives $i\geq D/2$, while $2i=d(w,w^{\prime})\leq D$ gives $i\leq D/2$. So $i=D/2$, which is excluded.
Assume now that $2i=\left\vert j-j^{\prime}\right\vert $.
Then $2i\leq D-(D-i)=i$, a contradiction.

Hence, at least one triple from $\{u,w,w^{\prime}\}$ and $\{v,w,w^{\prime}\}$ is non-collinear.
Consequently, there are at least $\binom{m_i}{2}$ triples with two vertices in $M_i(u)$ and the third being either $u$ or $v$.

\medskip

\noindent {\bf Case 2}:
{\em Triples containing $w,w^{\prime}\in M_D(u)$ and $u$.}

\smallskip

Since $2D>D$, these triples are non-collinear by Lemma \ref{Lemma_basic}$(ii)$.
There are $\binom{m_{D}}{2}$ triples of this type.

\medskip

\noindent {\bf Case 3}:
{\em (Only for even $D$.) Triples containing $w,w^{\prime}\in M_{D/2}(u)$. The third vertex is $u$, $v$ or it is also in $M_{D/2}(u)$.}

\smallskip

Take $w\not = w^{\prime}\in M_{D/2}(u)$ and put $r=d(w,w^{\prime})$, $a=d(v,w)$ and $b=d(v,w^{\prime})$.
Here $d(u,w)=D/2$ and $d(u,v)=D$, so the triangle inequality gives $a,b\geq D/2$ and $a,b\leq D$.
Hence $a+b\geq D$, with equality if and only if $a=b=D/2$, that is, if and only if $w,w^{\prime}\in M$.
Also $\left\vert a-b\right\vert \leq D/2<D$.

By the apex criterion at apex $u$, $\{u,w,w^{\prime}\}$ is collinear if and only if $r=D$.
By the apex criterion at apex $v$, $\{v,w,w^{\prime}\}$ is collinear if and only if $r=a+b$ or $r=\left\vert
a-b\right\vert $. Now $r\leq D\leq a+b$, so the first of these two options requires $r=D$
\emph{and} $w,w^{\prime}\in M$, while the second requires $r\leq D/2$.

This yields four families of non-collinear triples.
The first one consists of the triples $\{u,w,w^{\prime}\}$ for the $\binom{m^{\ast}}{2}-\left\vert E(H)\right\vert $ pairs with $r\not =D$.
The second one consists of the triples $\{v,w,w^{\prime}\}$ for the $\left\vert E(H)\right\vert -\left\vert E(J)\right\vert $ pairs with $r=D$ which do not lie both in $M$.
The third one consists of the triples $\{v,w,w^{\prime}\}$ for the $\binom{\left\vert M\right\vert}{2}-\left\vert E(J)\right\vert $ pairs inside $M$ with $r\not = D$.
The fourth one consists of the $t(H)$ triangles of $H$, which are triples of vertices at mutual distances $D$ and hence not collinear.
Obviously, the four families are pairwise disjoint, and the number of triples in these families is
$$
\Theta=\binom{m^{\ast}}{2}+\binom{\left\vert M\right\vert }{2}-2\left\vert E(J)\right\vert
+t(H).
$$

We now bound $\Theta$ from below.
Lemma~{\ref{Lemma_triangles}} applied to $J$, together with $t(H)\geq t(J)$ and $\left\vert M\right\vert \leq m^{\ast}$, gives $\Theta\geq \binom{m^{\ast}}{2}-(m^{\ast}+1)/2$.
And summing the bounds from all three cases we get the expression of the lemma.
\end{proof}

Now we bound the non-collinear triples of type $\tau_3$ and $\tau_4$.
Recall that a pair $w,w^{\prime}$ from a common level $M_i(u)$ is $u$-good if $d(w,w^{\prime})\not = 2i$.
By Lemma~{\ref{Lemma_basic}}(ii), if $2i>D$ then the pair $w,w^{\prime}$ is $u$-good and the triple $\{u,w,w^{\prime}\}$ is non-collinear.
In the following we abbreviate the notion $u$-good to good.

We consider vertex-disjoint good pairs $P_1,\dots,P_k$.
For a good pair $P_{j}=\{w_{j},w_{j}^{\prime}\}$ contained in a level $M_{i_{j}}(u)$ and for
$r_{j}=d(w_{j},w_{j}^{\prime})$ we denote by
$$
Z_{j}=V(G)\setminus (\{u,v\}\cup M_{i_{j}}(u)),\qquad E_{j}=\{z:d(u,z)\in \{r_{j}/2,\ i_{j}%
-r_{j}/2,\ i_{j}+r_{j}/2\}\}
$$
the set of vertices available to $P_{j}$ and the exceptional set of $P_{j},$ respectively.
By $L_{3}$ we denote the sum of the three largest values from $m_{1},\ldots ,m_{D}$.

\begin{lemma}
\label{Lemma_pairs}Let $G$ be a $2$-connected graph of diameter $D\geq 4$, let $(u,v)$ be a diametral pair of $G$, and let $P_{1},\ldots ,P_{k}$ be pairwise disjoint $u$-good pairs, each of them contained in a level $M_{i_{j}}(u)$ and avoiding $v$.
Then $\left\vert E_{j}\right\vert \leq L_{3}$ for every $j$, and the number of non-collinear triples of type $\tau_3$ and $\tau_4$ is at least
$$
\frac{1}{2}\sum_{j=1}^{k}\left( \left\vert Z_{j}\right\vert-\left\vert E_{j}\right\vert \right) .
$$
\end{lemma}

\begin{proof}
Since $r_{j}\geq 1$ and $r_{j}\not =2i_{j}$ by $u$-goodness, none of the numbers $r_{j}/2$, $i_{j}-r_{j}/2$, $i_{j}+r_{j}/2$ equals $0$.
Hence $E_{j}$ is contained in the union of at most three of the levels $M_{1}(u),\ldots ,M_{D}(u)$, and so $\left\vert E_{j}\right\vert \leq L_{3}$.

\medskip

\noindent\textbf{Claim A.} \emph{Let $w\not =w^{\prime}$ be two vertices of $M_{i}(u)$ with $r=d(w,w^{\prime})$, and let $z\not \in \{u,w,w^{\prime}\}$ with $b=d(u,z)$.
If $r\not \in \{2i,2b,2\left\vert b-i\right\vert \}$, then at least one of the triples $\{u,z,w\}$, $\{u,z,w^{\prime}\}$, $\{z,w,w^{\prime}\}$ is not collinear.}\smallskip

\noindent Assume that the first two triples are collinear.
By the apex criterion with apex $u$, each of $d(z,w),d(z,w^{\prime})$ lies in $\{b+i,\left\vert b-i\right\vert \}$.
By the triangle inequality $r\leq 2i$.
Assume that both these distances equal $\left\vert b-i\right\vert $.
Then the apex $z$ is equidistant from $w$ and $w^{\prime}$, so $\{z,w,w^{\prime}\}$ is collinear only if $r=2\left\vert b-i\right\vert $ or $r=0$.
Both are excluded. Assume next that both distances equal $b+i$.
Then collinearity needs $r=2(b+i)$ or $r=0$, while $r\leq 2i<2(b+i)$.
In the remaining case the two distances are $b+i$ and $\left\vert b-i\right\vert $, and they are distinct because $b,i\geq 1$.
Their sum is $2\max(b,i)$ and their difference is $2\min(b,i)$, so collinearity needs $r\in \{2b,2i\}$.
This is excluded as well, which establishes Claim A.

\medskip

\noindent Now fix $j$ and split $Z_{j}\setminus E_j$ into two parts.
Let $A_{j}$ be the set of those $z\in Z_j\setminus E_j$ for which at least one of $\{u,z,w_{j}\}$, $\{u,z,w_{j}^{\prime}\}$ is not collinear, and put $B_{j}=(Z_{j}\setminus E_j)\setminus A_{j}$.
Take first $z\in B_{j}$ and denote $b=d(u,z)$.
Both the above triples are then collinear, and $r_{j}\not \in \{2i_{j},2b,2\left\vert b-i_{j}\right\vert \}$.
So Claim A makes the triple $\{z\}\cup P_{j}$ non-collinear.
This triple avoids $u$, and it also avoids $v$ because $z\not =v$ and $v\not \in P_{j}$.
So this triple is of type $\tau_3$.
If $z\in A_{j}$ then there is a non-collinear triple of type $\tau_4$ since $Z_j$ does not contain vertices of $M_{i_{j}}(u)$.

Thus, for every $z\in Z_j\setminus E_j$ we get at least one non-collinear triple.
A triple of type $\tau_3$ is produced only once, since it contains only one of the pairs $P_{1},\ldots ,P_{k}$ (they are disjoint) and it does not contain $u$.
However, a triple $\{u,w,z\}$ of type $\tau_4$ (analogously $\{u,w^{\prime},z\}$) can be produced twice.
It is produced once by the pair containing $w$, and once by the pair containing $z$ (if such a pair exists).
The pairs are pairwise disjoint, so each of the two pairs is determined uniquely.
Hence every triple is produced at most twice.
This gives the lower bound $\frac 12\sum_{j=1}^{k}\left(\left\vert Z_{j}\right\vert -\left\vert E_{j}\right\vert \right)$.
\end{proof}

We are now in a position to prove the main estimate of this section.

\begin{theorem}
\label{Tm_2conn}
Every $2$-connected graph on $n\geq 31$ vertices has at least $2(n-4)$ non-collinear triples.
\end{theorem}

\begin{proof}
For $D\leq 3$ the claim is Theorem \ref{Tm_diam3}.
So assume that $D\geq 4$, and fix a diametral pair $(u,v)$.
We use the following abbreviations
$$
m=n-4,\qquad s=D-3,\qquad p=\left\lceil D/2\right\rceil -1,\qquad R=n-1-L_{3}.
$$
Notice that $s\geq 1$.
Every vertex other than $u$ lies in exactly one of the $D$ levels $M_{i}(u)$, $1\leq i\leq D$.
Hence $m_{1}+\cdots +m_{D}=n-1=m+3$.
Each of the $D$ levels is non-empty, so $R\geq s\geq 1$.
Here $R$ is the total size of the $s$ levels which are not among the three largest ones, and $L_{3}=m+3-R$.

\medskip

\noindent\textbf{Step 1.}
We call a level $M_{i}(u)$ \emph{far} if $D/2<i<D$.
There are exactly $p$ far levels.
Each of them has at least two vertices by Lemma \ref{Lemma_basic}$(iii)$.
In every far level we choose one pair of vertices.
The resulting $p$ pairs are pairwise disjoint, and they avoid $v$, since $v\in M_{D}(u)$.
They are $u$-good by Lemma \ref{Lemma_basic}$(ii)$.
Let $P_{j}$ be the pair chosen in the far level $M_{i}(u)$, and write $i=i_{j}$.
Since $u,v\not \in M_{i}(u)$, we have $\left\vert Z_{j}\right\vert =n-2-m_{i}$.
Lemma \ref{Lemma_pairs} gives $\left\vert E_{j}\right\vert \leq L_{3}$, and so
$\left\vert Z_{j}\right\vert -\left\vert E_{j}\right\vert \geq n-2-m_{i}-L_{3}=R-1-m_{i}$.
Hence, by Lemma \ref{Lemma_pairs}, the number of non-collinear triples of type $\tau_{3}$ and $\tau_{4}$ is at least
$$
\frac{1}{2}\sum_{D/2<i<D}\left( R-1-m_{i}\right) =\frac{p(R-1)}{2}-\frac{1}{2}\sum_{D/2<i<D}m_{i}.
$$

\medskip

\noindent\textbf{Step 2.}
The triples counted in Lemma \ref{Lemma_levels} and those counted in Lemma \ref{Lemma_pairs} are disjoint.
So we may add the bound of Lemma \ref{Lemma_levels} to the bound of Step 1.
We use $\binom{x}{2}=\frac{x(x-2)}{2}+\frac{x}{2}$ for every level.
Let $N$ be the set of indices $i\in \{1,\ldots ,D\}$ such that $M_{i}(u)$ is neither far nor the middle level $M_{D/2}(u)$.
Then $\frac{1}{2}\sum_{i=1}^{D}m_{i}-\frac{1}{2}\sum_{D/2<i<D}m_{i}=\frac{1}{2}\sum_{i\in N}m_{i}+\frac{m^{\ast }}{2}$, and it is true also if $D$ is odd, since then $m^{\ast }=0$.
We obtain
$$
\mathrm{noncol}(G)\geq \sum_{i=1}^{D}\frac{m_{i}(m_{i}-2)}{2}+\frac{1}{2}\sum_{i\in N}m_{i}-\frac{1}{2}+\frac{p(R-1)}{2}.
$$
Notice that $N=\{1,\ldots ,p\}\cup \{D\}$.
By Lemma \ref{Lemma_basic}$(iii)$ we have $m_{i}\geq 2$ for $1\leq i\leq p$, and $m_{D}\geq 1$ since $v\in M_{D}(u)$.
Hence $\sum_{i\in N}m_{i}\geq 2p+1$.
Since also $\sum_{i=1}^{D}m_{i}(m_{i}-2)/2=\frac{1}{2}\sum_{i=1}^{D}m_{i}^{2}-(m+3)$, we get
\begin{equation}
\mathrm{noncol}(G)\geq \frac{1}{2}\sum_{i=1}^{D}m_{i}^{2}-(m+3)+\frac{p(R+1)}{2}.
\label{For_2connStar}
\end{equation}

\medskip

\noindent\textbf{Step 3.}
We split the $D$ levels into the three largest ones, of total size $L_{3}$, and the remaining $s$ ones, of total size $R$.
The Cauchy--Schwarz inequality $\sum_{i=1}^{k}a_{i}^{2}\geq (a_{1}+\cdots +a_{k})^{2}/k$, applied to each of the two groups, gives
$\sum_{i=1}^{D}m_{i}^{2}\geq L_{3}^{2}/3+R^{2}/s$.
Next, $p\geq (s+1)/2$.
Indeed, if $D$ is even then $p=D/2-1=(s+1)/2$, and if $D$ is odd then $p=(D-1)/2=(s+2)/2$.
Since $R+1>0$, inequality (\ref{For_2connStar}) yields
$$
\mathrm{noncol}(G)\geq \frac{L_{3}^{2}}{6}+\left( \frac{R^{2}}{2s}+\frac{s(R+1)}{4}\right) +\frac{R+1}{4}-(m+3).
$$

\medskip

\noindent\textbf{Step 4.}
By the AM--GM inequality applied to the two summands in the bracket, the bracket is at least
$2\sqrt{R^{2}(R+1)/8}=R\sqrt{(R+1)/2}$.
We show that
$$
R\sqrt{(R+1)/2}\geq \frac{53R-289}{12}.
$$
If $53R\leq 289$, then the right-hand side is not positive, and the inequality holds.
Otherwise both sides are positive.
Squaring and multiplying by $144$, the inequality becomes $72R^{2}(R+1)\geq (53R-289)^{2}$.
Notice that $72R^{2}(R+1)-(53R-289)^{2}=(R-17)^{2}(72R-289)$.
This is non-negative, since $72R>53R>289$.
Hence
$$
\mathrm{noncol}(G)\geq \frac{L_{3}^{2}}{6}+\frac{53R-289}{12}+\frac{R+1}{4}-(m+3)=\frac{L_{3}^{2}}{6}+\frac{28R-143}{6}-(m+3).
$$
We substitute $R=m+3-L_{3}$.
The right-hand side becomes $\left( L_{3}^{2}-28L_{3}+22m-77\right) /6=\left( (L_{3}-14)^{2}+22m-273\right) /6$.
Consequently
$$
\mathrm{noncol}(G)-2m\geq \frac{(L_{3}-14)^{2}+10m-273}{6}\geq \frac{10m-273}{6},
$$
and for $m\geq 27$ the last expression is at least $-\frac 12$.
Since $\mathrm{noncol}(G)-2m$ is an integer, it must be non-negative.
So $\mathrm{noncol}(G)\geq 2m=2(n-4)$ whenever $n\geq 31$.
\end{proof}

Combining Proposition~{\ref{Prop_dumbbell}}, Theorem~{\ref{Tm_cutvertex}} and Theorem~{\ref{Tm_2conn}} we obtain the next statement.

\begin{corollary}
\label{Cor_pendant}Among all connected graphs $G$ on $n\geq 31$ vertices without pendant vertices it
holds that
\[
\mathrm{Int}(G)\leq \frac{n^{3}-n}{6}-2(n-4),
\]
and equality is attained by the dumbbell $D_{n}.$
\end{corollary}

\begin{proof}
Notice that $\binom{n}{2}+\binom{n}{3}=(n^{3}-n)/6$, so by Lemma \ref{Lemma_identity} we have $\mathrm{Int}(G)=(n^{3}-n)/6-\mathrm{noncol}(G)$.
A connected graph on $n\geq 3$ vertices without pendant vertices either has a cut vertex, or it is $2$-connected.
In both cases $\mathrm{noncol}(G)\geq 2(n-4)$ by Theorem \ref{Tm_cutvertex} and Theorem \ref{Tm_2conn}, respectively.
By Proposition \ref{Prop_dumbbell}, $\mathrm{noncol}(D_{n})=2(n-4)$.
\end{proof}

\section{Concluding remarks}

We have shown that
$$
\mathrm{Int}(G)\leq \frac{n^{3}-n}{6}
$$
for every connected graph $G$ on $n$ vertices.
Moreover, equality forces $G\cong P_{n}$ with the single exception $G\cong C_{4}$ at $n=4$, where $\mathrm{Int}(C_4)=\mathrm{Int}(P_4)$.

As regards Problem~{\ref{Prob_pendant}}, we have shown that among all graphs $G$ on
$n\geq 31$ vertices without pendant vertices it holds that
$$
\mathrm{Int}(G)\leq \frac{n^{3}-n}{6}-2(n-4),
$$
and equality is attained by the dumbbell graph consisting of two copies of $C_4$ joined by a path.

Finally, the identity of Lemma \ref{Lemma_identity} reduces every extremal question about $\mathrm{Int}$ to a question about collinear triples.
It would be interesting to see which further extremal problems for $\mathrm{Int}$ become accessible in this way.
In particular, the remaining problem of \cite{IntervalIndex}, asking which graphs on $2m$ vertices satisfy $\mathrm{Int}(G)=\mathrm{Int}(C_{2m})$, is by Lemma~{\ref{Lemma_identity}} nothing else but the question which graphs on $2m$ vertices have exactly $2\binom{m}{3}$ non-collinear triples.

\bigskip

\bigskip\noindent\textbf{Declaration of AI use.}~~The authors used Claude (Anthropic) to check and simplify the proofs and to improve the presentation.
The authors reviewed all such content and take full responsibility for the paper.

\emergencystretch=1.5em
\bigskip\noindent\textbf{Acknowledgments.}~~M.~Knor acknowledges partial support by Slovak
research grants APVV-23-0076, APVV-22-0005, VEGA 1/0069/23 and VEGA 1/0011/25.
R.~\v{S}krekovski acknowledges partial support by ARIS projects J1-3002 and J1-70016, program P1-0383, and the annual work program of Rudolfovo, and by the bilateral Croatian--\allowbreak Slovenian project BI-HR/25-27-004.
J.~Sedlar acknowledges partial support by the bilateral Croatian--\allowbreak Slovenian project BI-HR/25-27-004, by the Croatian Science Foundation project HRZZ-IP-2024-05-2130, and by the NextGenerationEU instrument via IP-UNIST-17 (GEORAZ).


\begin{thebibliography}{99}    %

\bibitem {DA}A.~Alochukwu, P.~Dankelmann, Wiener index in graphs with given minimum degree and
maximum degree, \emph{Discrete Math.\ Theor.\ Comput.\ Sci.}\ \textbf{23} (2021) \#11.

\bibitem {BDHKS}S.~Bessy, F.~Dross, K.~Hri\v{n}\'{a}kov\'{a}, M.~Knor, R.~\v{S}krekovski,
Maximal Wiener index for graphs with prescribed number of blocks,
\emph{Appl.\ Math.\ Comput.}\ \textbf{380} (2020) 125274.

\bibitem {BDKS}S.~Bessy, F.~Dross, M.~Knor, R.~\v{S}krekovski, Graphs with the second and third
maximum Wiener indices over the $2$-vertex connected graphs, \emph{Discrete Appl.\ Math.}\
\textbf{284} (2020) 195--200.

\bibitem {IntervalIndex}M.~Knor, \v{Z}.~Kovijani\'{c} Vuki\'{c}evi\'{c}, R.~\v{S}krekovski,
V.~Do\v{s}ljak, D.~Dimitrov, G.~Popivoda, S.~Vujo\v{s}evi\'{c}, J.~Sedlar, Interval index
of a graph, \emph{Discuss.\ Math.\ Graph Theory} (2026), in press, doi:10.7151/dmgt.2654.

\bibitem {KSY1}M.~Knor, J.~Sedlar, R.~\v{S}krekovski, Y.~Yang, Invitation to the subpath
number, \emph{Appl.\ Math.\ Comput.}\ \textbf{509} (2026) 129646.

\bibitem {KSY2}M.~Knor, J.~Sedlar, R.~\v{S}krekovski, Y.~Yang, The subpath number of cactus
graphs, \emph{Comput.\ Appl.\ Math.}\ \textbf{45} (2026) 102.

\bibitem {KSZ}M.~Knor, J.~Sedlar, R.~\v{S}krekovski, X.-D.~Zhang, Counting geodesic paths in
graphs, \emph{Mediterr.\ J.\ Math.}\ \textbf{23} (2026) 171.

\bibitem {KST}M.~Knor, R.~\v{S}krekovski, A.~Tepeh, Mathematical aspects of Wiener index,
\emph{Ars Math.\ Contemp.}\ \textbf{11} (2016) 327--352.

\bibitem {MoonMoser}J.~W.~Moon, L.~Moser, On a problem of Tur\'{a}n, \emph{Magyar Tud.\ Akad.\
Mat.\ Kutat\'{o} Int.\ K\"{o}zl.}\ \textbf{7} (1962) 283--286.

\bibitem {Mulder}H.~M.~Mulder, \emph{The Interval Function of a Graph}, Mathematical Centre
Tracts \textbf{132}, Mathematisch Centrum, Amsterdam, 1980.

\bibitem {Nebesky}L.~Nebesk\'{y}, The interval function of a connected graph and a
characterization of geodetic graphs, \emph{Math.\ Bohem.}\ \textbf{126} (2001) 247--254.

\bibitem {Wiener}H.~Wiener, Structural determination of paraffin boiling points, \emph{J.\ Amer.\
Chem.\ Soc.}\ \textbf{69} (1947) 17--20.
\end{thebibliography}
\end{document}